\documentclass[12pt,a4paper]{amsart}

\usepackage{euscript,amsfonts,amssymb,amsmath,amscd,amsthm,enumerate,hyperref}

\usepackage{comment}

\usepackage{mathrsfs}

\usepackage{graphicx}

\usepackage{color}

\usepackage{tikz}
\usepackage{pgf}
\usetikzlibrary{matrix,arrows,decorations.markings}
\input xy
\xyoption{all}

\usepackage[margin=1.1in]{geometry}

\usepackage{bbm}

\newtheorem{theorem}{Theorem}[section]
\newtheorem*{theorem*}{Theorem}
\newtheorem{lemma}[theorem]{Lemma}
\newtheorem{definition}[theorem]{Definition}
\newtheorem{proposition}[theorem]{Proposition}

\theoremstyle{remark}
\newtheorem{rmk}[theorem]{Remark}

\newcommand{\E}{\mathbb E}

\newcommand{\rr}{\mathbb{R}}

\newcommand{\GN}{\mathcal{G}^{(N)}}
\newcommand{\ep}{\hfill \ensuremath{\Box}}
\newcommand{\eq}{\begin{equation}}
\newcommand{\en}{\end{equation}}

\numberwithin{equation}{section}

\title[Dyson Brownian motions]{Intertwinings of beta-Dyson Brownian motions of different dimensions}

\author{Kavita Ramanan} 
\address{Division of Applied Mathematics, Brown University, Providence, RI, USA}
\email{kavita\_ramanan@brown.edu}

\author{Mykhaylo Shkolnikov}
\address{ORFE Department, Princeton University, Princeton, NJ, USA}
\email{mshkolni@gmail.com}

\begin{document}

\begin{abstract}
We show that for all positive $\beta$ the semigroups of $\beta$-Dyson Brownian motions of different dimensions are intertwined. The proof relates $\beta$-Dyson Brownian motions directly to Jack symmetric polynomials and omits an approximation of the former by discrete space Markov chains, thereby disposing of the technical assumption $\beta\ge1$ in \cite{GS}. The corresponding results for $\beta$-Dyson Ornstein-Uhlenbeck processes are also presented.
\end{abstract}

\thanks{The first author, K. Ramanan, is partially supported by NSF grant DMS-1407504, and the second author, M. Shkolnikov, is partially supported by NSF grant DMS-1506290.} 



\maketitle

\section{Introduction}

For a fixed $\beta>0$ consider the system of stochastic differential equations
\begin{equation}\label{eq:DBM}
\mathrm{d}X^{(N)}_i(t)=\frac{\beta}{2}\,\sum_{\substack{1\le j\le
    N\\j\neq i}}
\frac{1}{X^{(N)}_i(t)-X^{(N)}_j(t)}\,\mathrm{d}t+\mathrm{d}B^{(N)}_i(t),\quad
i=1,2,\ldots,N, 
\end{equation}
with initial condition satisfying $X^{(N)}_1(0)\le
X^{(N)}_2(0)\le\cdots\le X^{(N)}_N(0)$ and where
$B^{(N)}=(B^{(N)}_1,B^{(N)}_2,\ldots,B^{(N)}_N)$ is a standard
$N$-dimensional Brownian motion. For $\beta=1,2,4$ it has been shown
by Dyson \cite{Dy} that the equation \eqref{eq:DBM} describes the
evolution of eigenvalues of random symmetric matrices whose entries
follow real, complex, and quaternion Brownian motions,
respectively. More generally, one can make sense of the unique strong
solution to \eqref{eq:DBM} for all $\beta>0$ (see \cite[Theorem
3.1]{CL1} and also the alternative proof in \cite[Proposition
4.3.5]{AGZ} for $\beta\ge1$). We refer to the latter as the
$N$-dimensional $\beta$-Dyson Brownian motion and write $P^{(N)}(t)$,
$t\ge0,$ for the associated Markov transition operators, emphasizing
the dependence on dimension and omitting the explicit dependence on $\beta$. 

\medskip

More recently, Dyson Brownian motions have been used to establish a variety of universality conjectures in random matrix theory, such as  universality of the local eigenvalue statistics in the bulk and at the edge of the spectrum of symmetric Wigner matrices (see \cite{ESY} and \cite{BEY}, respectively),  universality of the local statistics in the bulk and at the edge of non-critical beta ensembles with smooth external potentials (see \cite{BEY}) and the universality of the bulk statistics and the distribution of the second largest eigenvalue of dense Erd\"os-R\'{e}nyi random graphs (see \cite{EKYY1}, \cite{EKYY2}) among others.   

\medskip

Concurrently, it was observed in \cite{Ba} that for $\beta=2$ and
densely packed initial conditions
$X^{(N)}_1(0)=X^{(N)}_2(0)=\cdots=X^{(N)}_N(0),$ the largest
coordinate process
$X^{(N)}_N$ in an $N$-dimensional Dyson Brownian motion has the same distribution as the rightmost particle in the Brownian totally asymmetric simple exclusion process (TASEP) started from the same densely packed initial condition. This connection reveals that the fixed time fluctuations of the rightmost particle in the Brownian TASEP are governed by the same distribution as those of the largest eigenvalue of a Hermitian random matrix, namely the Tracy-Widom distribution $TW_2$. 

\medskip

The seminal paper \cite{Wa} explains the relation between the $(\beta=2)$-Dyson Brownian motion and the Brownian TASEP by constructing a process $(X^{(k)}_i,\,1\le i\le k\le N)$ taking values in the Gelfand-Tseitlin cone 
\begin{equation}
\GN:=\Big\{\big(x^{(k)}_i:1\le i\le k\le N\big)\in\rr^{N(N+1)/2}:\,x^k_i\le x^{k-1}_i\le x^k_{i+1},\,1\le i\le k\le N-1\Big\}
\end{equation}
and starting from $0\in\GN$ such that each ``level'' $X^{(k)}:=(X^{(k)}_1,X^{(k)}_2,\ldots,X^{(k)}_k)$ performs a $k$-dimensional $(\beta=2)$-Dyson Brownian motion, whereas the ``diagonal'' $(X^{(1)}_1,X^{(2)}_2,\ldots,X^{(N)}_N)$ forms the Brownian TASEP with $N$ particles. We refer to $(X^{(k)}_i,\,1\le i\le k\le N)$ as the $(\beta=2)$-multilevel Dyson Brownian motion. The construction in \cite{Wa} relies on the fact that for $\beta=2$ and any $k\ge1$ the semigroups $P^{(k)}(t)$, $t\ge0,$ and $P^{(k+1)}(t)$, $t\ge0,$ are \textit{intertwined}. More specifically, there exists a stochastic kernel $L^{(k)}(x^{(k+1)},\mathrm{d}x^{(k)})$ such that 
\begin{equation}\label{Warren_intertw}
L^{(k)}\,P^{(k)}(t)=P^{(k+1)}(t)\,L^{(k)},\quad t\ge0. 
\end{equation}
Figure \ref{fig:intertwining} provides an illustration of the intertwining relationship \eqref{Warren_intertw}.

\begin{figure}[t]
\centerline{
\xymatrix@=6em{
X^{(k+1)}(0) \ar[d]_{L^{(k)}} \ar[r]^{P^{(k+1)}(t)} & X^{(k+1)}(t) \ar[d]^{L^{(k)}} \\
X^{(k)}(0) \ar[r]^{P^{(k)}(t)} & X^{(k)}(t)}
} 
\caption{Intertwining of $P^{(k)}$ and $P^{(k+1)}$.}
\label{fig:intertwining}
\end{figure} 

\medskip

By now intertwinings have been used to construct a number of important ``multilevel'' processes. For discrete space Markov chains the procedure is based on the coupling of Diaconis and Fill \cite{DF} and has been applied to construct multilevel Markov chains describing the growth of two-dimensional random surfaces, see \cite{BC}, \cite{BF2}, \cite{BK}, \cite{GS} and the references therein. Such random surfaces are of great interest as members of the $(2+1)$-dimensional Kardar-Parisi-Zhang (KPZ) universality class and yield insights into solutions of the KPZ stochastic partial differential equation.   

\medskip

The construction of multilevel processes with continuous state space presents additional technical challenges, because the transition operators of the underlying Markov processes are usually not available in closed form (the situation of \cite{Wa} being a notable exception). For this reason, there are only three known examples in the continuous space setting: the two-dimensional Whittaker growth model of \cite{OC}, the $\beta=2$ multilevel Dyson Brownian motion of \cite{Wa} and the $\beta>2$ multilevel Dyson Brownian motions of \cite{GS}. The article \cite{OC} is the only one where the intertwining relationship behind the multilevel continuous space process is established directly and in the absence of explicit formulas for the transition operators. In contrast, \cite{Wa} relies on such explicit formulas and the proof in \cite{GS} is based on an approximation by discrete space Markov chains.   

\medskip

In this article we make progress towards the construction of
multilevel Dyson Brownian motions for all positive $\beta$ by
extending the intertwining relationship \eqref{Warren_intertw} to all
$\beta>0$.   In contrast with  \cite{GS}, our proof works directly with $\beta$-Dyson Brownian motions rather than their discrete space approximations. In the process, we provide formulas for how $\beta$-Dyson Brownian motion generators and semigroups act on Jack symmetric polynomials. 

\medskip

To state our main result we need the following set of notations. We 
define the chambers
\begin{eqnarray}
&&\;\mathcal{W}^{(k+1)}:=\big\{x^{(k+1)}\in\rr^{k+1}:\;x^{(k+1)}_1\le x^{(k+1)}_2\le\cdots\le x^{(k+1)}_{k+1}\big\}, \\
&&\;\mathcal{W}^{(k)}(x^{(k+1)}):=\big\{x^{(k)}\in\rr^k:
x^{(k+1)}_1\le x^{(k)}_1\le x^{(k+1)}_2\le\cdots\le x^{(k)}_k\le
x^{(k+1)}_{k+1}\big\}, 
\end{eqnarray}
for $x^{(k+1)}\in\mathcal{W}^{(k+1)}$. 
In addition, we fix $\beta > 0$ and for each $x^{(k+1)}\in\mathcal{W}^{(k+1)}$ we write 
\begin{equation}\label{DA_probdens}
\begin{split}
\Lambda^{(k)}(x^{(k+1)},x^{(k)}):=\frac{\Gamma(\beta(k+1)/2)}{\Gamma(\beta/2)^{k+1}}\,\prod_{1\le
  i<j\le k+1} \big(x^{(k+1)}_j-x^{(k+1)}_i\big)^{1-\beta}\,\prod_{1\le
  i<j\le k} \\
\big(x^{(k)}_j-x^{(k)}_i\big) \prod_{i=1}^k \prod_{j=1}^{k+1} \big|x^{(k)}_i-x^{(k+1)}_j\big|^{\beta/2-1}
\end{split}
\end{equation}
for the Dixon-Anderson conditional probability density of $x^{(k)}$ on
$\mathcal{W}^{(k)}(x^{(k+1)})$ given $x^{(k+1)}$, where $\Gamma(\cdot)$ stands for the gamma function. The density $\Lambda^{(k)}(x^{(k+1)},\,\cdot\,)$ has been introduced independently in \cite{Di} and \cite{An} and for $\beta=1,2,4$ describes the conditional distribution of the eigenvalues of the $k\times k$ corner in a $(k+1)\times(k+1)$ random matrix from the Gaussian orthogonal, unitary, symplectic ensemble, respectively (see also \cite{Fo1} for a more detailed discussion and generalizations). 

\medskip

Our main result reads as follows. 

\begin{theorem}\label{main_thm}
Let $\beta>0$. Then, with the Dixon-Anderson conditional probability
density of \eqref{DA_probdens}, the stochastic kernel
$L^{(k)}(x^{(k+1)},\mathrm{d}x^{(k)}):=\Lambda^{(k)}(x^{(k+1)},x^{(k)})\,\mathrm{d}x^{(k)}$
intertwines the semigroups $P^{(k)}$ and $P^{(k+1)}$ of the
$k$-dimensional and the $(k+1)$-dimensional $\beta$-Dyson Brownian
motions, respectively.  In other words, one has the following
equalities of probability measures: 
\begin{equation}\label{main_intertw}
L^{(k)}\,P^{(k)}(t)=P^{(k+1)}(t)\,L^{(k)},\quad t\ge0. 
\end{equation}
\end{theorem}

\begin{rmk}\label{rmk:iter}
An iteration of \eqref{main_intertw} shows that for all $1\le k<n$,
the following equality holds  
\begin{equation}
\label{prod-intertw}
\left( \prod_{m=k}^{n-1} L^{(m)}\right)\,P^{(k)}(t)=P^{(n)}(t)\,\left(\prod_{m=k}^{n-1} L^{(m)}\right),\quad t\ge0.
\end{equation}
\end{rmk}

\begin{rmk}\label{rmk:coupling}
The general results on intertwinings of diffusion processes in \cite{PS} (see, in particular, \cite[Theorem 5]{PS}) suggest that one should be able to realize the intertwining relationships \eqref{main_intertw} for $k=1,2,\ldots,N-1$ by a coupling of $\beta$-Dyson Brownian motions of dimensions $1,2,\ldots,N$ to a multilevel process on $\GN$. Moreover, an iterative application of the formulas in \cite[Theorem 5]{PS} leads one to conjecture that the generator of such a coupling should be given by 
\begin{equation}\label{multi_gen}
\frac{1}{2}\,\sum_{1\le i\le k\le N} \partial_{x^{(k)}_i}^2 
-\sum_{1\le i\le k\le N} \sum_{j\neq i} \frac{\beta/2-1}{x^{(k)}_i-x^{(k)}_j}\,\partial_{x^{(k)}_i}
+\sum_{1\le i\le k+1\le N} \sum_{j=1}^k \frac{\beta/2-1}{x^{(k+1)}_i-x^{(k)}_j}\,\partial_{x^{(k+1)}_i},
\end{equation}
endowed with the reflecting boundary conditions described in
\cite{Wa}. Processes of this type have indeed been constructed in
\cite{Wa} for $\beta=2$ and in \cite{GS} for $\beta>2$ and one should
view them as $\GN$-valued analogues of Bessel processes of dimensions
$d=1$ and $d>1$, respectively. This analogy suggests that for
$\beta<2$ a process with generator \eqref{multi_gen} should no longer
be a semimartingale.  The construction of such a process appears
challenging: the non-reversibility of the process in conjunction with
the reflecting boundary conditions rule out non-symmetric Dirichlet
form techniques (note that the process constructed in \cite[Exercise
II.2.14, Theorem IV.3.5]{MR} is absorbed at the boundary); on the
other hand, the singularity of the coefficients together with the
inapplicability of localization techniques as in \cite{GS} prevent one
from directly building the coupling from known processes. 
\end{rmk}

\smallskip

We conclude the introduction by presenting the version of Theorem \ref{main_thm} for $\beta$-Dyson Ornstein-Uhlenbeck processes. Fix a $\beta>0$ and consider the system of stochastic differential equations 
\begin{equation}\label{eq:DOU}
\mathrm{d}Y^{(N)}_i(t)=\frac{\beta}{2}\,\sum_{\substack{1\le j\le N\\j\neq i}} \frac{1}{Y^{(N)}_i(t)-Y^{(N)}_j(t)}\,\mathrm{d}t-\frac{Y^{(N)}_i(t)}{2}\,\mathrm{d}t+\mathrm{d}B^{(N)}_i(t),
\end{equation}
$i=1,2,\ldots,N$,  with initial condition $Y^{(N)}_1(0)\le Y^{(N)}_2(0)\le\cdots\le Y^{(N)}_N(0)$ and where $B^{(N)}=(B^{(N)}_1,B^{(N)}_2,\ldots,B^{(N)}_N)$ is a standard $N$-dimensional Brownian motion as before. We call the unique strong solution of \eqref{eq:DOU} (see \cite[Theorem 3.1]{CL1}) the $N$-dimensional $\beta$-Dyson Ornstein-Uhlenbeck process and denote the associated Markov transition operators by $Q^{(N)}(t)$, $t\ge0$. 
  
\smallskip

\begin{proposition}\label{prop_OU}
Let $\beta>0$. Then, with the notations of Theorem \ref{main_thm} and the previous paragraph, one has the intertwining relationships
\begin{equation}\label{OU_intertw_semi}
L^{(k)}\,Q^{(k)}(t)=Q^{(k+1)}(t)\,L^{(k)},\quad t\ge0, 
\end{equation}
for all $k\ge1$. 
\end{proposition}
  
\begin{rmk}
As in Remark \ref{rmk:iter}, 
\eqref{prod-intertw} holds with $P^{(k)}$ and $P^{(n)}$ replaced by
$Q^{(k)}$ and $Q^{(n)}$, respectively. 
\end{rmk}

The rest of the paper is structured as follows. In Section
\ref{sec:Jacks} we present some facts about Jack symmetric polynomials
that are needed in the proof of Theorem \ref{main_thm}. Section \ref{sec:main} is then devoted to the proof of Theorem \ref{main_thm}. In Section \ref{sec:beta1} we present a simpler proof for the case $\beta=1$ based on the random matrix interpretation of $\beta=1$ Dyson Brownian motions. Finally, in Section \ref{sec:OU} we give the proof of Proposition \ref{prop_OU}.

\section{Preliminaries on Jack symmetric polynomials} \label{sec:Jacks}

We start with the definition of Jack symmetric polynomials following \cite[Section 2]{BF1}, but replacing the parameter $\alpha$ there by $\theta:=1/\alpha$ as in \cite{OO}. 

\begin{definition}\label{def:Jack}
Let $\theta:=\beta/2$ and consider the differential operator
\begin{equation}\label{op_Jack}
\sum_{i=1}^k z_i^2\,\frac{\partial^2}{\partial z_i^2} + 2\theta\,\sum_{i=1}^k \, \sum_{\substack{1\le j\le k\\j\neq i}} \frac{z_i^2}{z_i-z_j}\,\frac{\partial}{\partial z_i}
\end{equation}
acting on symmetric polynomials in $k$ variables. It is known (see
\cite[Theorem 3.1]{St}) that the eigenfunctions of this operator can
be indexed by non-decreasing positive integer sequences $\kappa =
(\kappa_i, i \leq \ell)$, with $\kappa_1\ge\kappa_2\ge\cdots\ge\kappa_l>0$ and $l\le k$, such that the lexicographically leading monomial with a non-zero coefficient in the eigenfunction $J_\kappa(\,\cdot\,;\theta)$ is $z_1^{\kappa_1}z_2^{\kappa_2}\cdots z_l^{\kappa_l}$ and $J_\kappa(\,\cdot\,;\theta)$ is normalized according to
\begin{equation}\label{norm_formula}
J_\kappa(1_k;\theta)=\prod_{i=1}^l \,\prod_{j=1}^{\kappa_i} \,(k+1-i+(j-1)/\theta)=\theta^{-|\kappa|}\,\prod_{i=1}^l \frac{\Gamma(\!(k+1-i)\theta+\kappa_i)}{\Gamma(\!(k+1-i)\theta)}.
\end{equation}
Here, $1_k$ is the $k$-dimensional vector whose components are all equal to $1$ and $|\kappa|:=\sum_{i=1}^l \kappa_i$. The eigenfunction $J_\kappa(\,\cdot\,;\theta)$ is called the \textit{Jack symmetric polynomial} in $k$ variables with parameters $\kappa$, $\theta$. 
\end{definition}

\begin{rmk}
Jack symmetric polynomials can be also defined as the eigenfunctions of the Sekiguchi differential operators
\begin{equation}\label{Sekiguchi}
\frac{1}{\prod_{1\le i<j\le k} (z_i-z_j)}\,\det\bigg[z_i^{k-j}\bigg(z_i\,\frac{\partial}{\partial z_i}+(k-j)\,\theta+u\bigg)\bigg]_{1\le i,j\le k}
\end{equation}
in the space of symmetric polynomials in $k$ variables (see, e.g., \cite[Section 1]{OO} and the references therein). The equivalence of the two definitions follows from the fact that the operator \eqref{op_Jack} can be recovered by an affine transformation of the operator multiplying $u^{k-2}$ in the expansion of \eqref{Sekiguchi} in powers of the auxiliary variable $u$ (see \cite[Example VI.3.3]{Ma}). 
\end{rmk}

\smallskip

We proceed to the definition of the generalized binomial coefficients associated with shifts of Jack symmetric polynomials (see \cite[remark on p.73]{OO} for more details).

\begin{definition}
The coefficients $\binom{\kappa}{\rho}_\theta$ in the expansion 
\begin{equation}\label{eq:binomial}
\frac{J_\kappa(1_k+z;\theta)}{J_\kappa(1_k;\theta)}=\sum_{m=0}^{|\kappa|}\,\sum_{|\rho|=m} \binom{\kappa}{\rho}_\theta\,\frac{J_\rho(z;\theta)}{J_\rho(1_k;\theta)}
\end{equation}
are referred to as the generalized binomial coefficients.
\end{definition}

\begin{rmk}
Note that the definition of $\binom{\kappa}{\rho}_\theta$ does not depend on the particular normalization of the Jack symmetric polynomials. 
\end{rmk}

\smallskip

Next, we state three identities for Jack symmetric polynomials from \cite[equations (2.13a), (2.13b), (2.13d)]{BF1} that will be employed in the proofs of Theorem \ref{main_thm} and Proposition \ref{prop_OU}.

\begin{proposition}\label{prop_B1B2}
The differential operators
\begin{equation}\label{B3def}
\begin{array}{rcl}
{\mathcal B}_1 & := & \displaystyle \sum_{i=1}^k \frac{\partial}{\partial z_i}, \\

{\mathcal B}_2 & := & \displaystyle\frac{1}{2}\,\sum_{i=1}^k
z_i\,\frac{\partial^2}{\partial z_i^2}+\theta\,\sum_{i=1}^k
\sum_{\substack{1\le j\le k\\j\neq i}}
\frac{z_i}{z_i-z_j}\,\frac{\partial}{\partial z_i}, \\
{\mathcal B}_3 & := & \displaystyle\sum_{i=1}^k z_i\,\frac{\partial}{\partial z_i}, 
\end{array}
\end{equation}
act on Jack symmetric polynomials in $k$ variables with parameter $\theta$ as follows:
\begin{eqnarray}
&& {\mathcal B}_1\,J_\kappa(z;\theta)=J_\kappa(1_k;\theta)\,\sum_{i=1}^l \binom{\kappa}{\kappa_{(i)}}_\theta\,\frac{J_{\kappa_{(i)}}(z;\theta)}{J_{\kappa_{(i)}}(1_k;\theta)}, \\
&& {\mathcal B}_2\,J_\kappa(z;\theta)
=J_\kappa(1_k;\theta)\,\sum_{i=1}^l \binom{\kappa}{\kappa_{(i)}}_\theta\, (\kappa_i-1+(k-i)\,\theta)\,\frac{J_{\kappa_{(i)}}(z;\theta)}{J_{\kappa_{(i)}}(1_k;\theta)}, \\
&& {\mathcal B}_3\,J_\kappa(z;\theta)=|\kappa|\,J_\kappa(z;\theta).
\label{B3action}
\end{eqnarray}
Here $\kappa_{(i)}$ is the sequence obtained from $\kappa$ by replacing $\kappa_i$ by $\kappa_i-1$ unless $i=l$ and $\kappa_l=1$ in which case we drop $\kappa_l$ from $\kappa$. We have also set $\binom{\kappa}{\kappa_{(i)}}_\theta=0$ whenever $\kappa_{(i)}$ is no longer a non-decreasing positive sequence. 
\end{proposition}

\smallskip

The last ingredient we need is a formula for the action of the Dixon-Anderson conditional probability density of \eqref{DA_probdens} on the Jack symmetric polynomials in $k$ variables. A proof of this formula can be found in \cite[Section 6]{OO} (note that the particular normalization of the Jack symmetric polynomials is irrelevant here). 

\begin{proposition}\label{prop_DA_OO}
With $\Lambda^{(k)}$ of \eqref{DA_probdens} and $\theta=\beta/2$ one
has for every $x^{(k+1)} \in {\mathcal W}^{(k+1)}$ and sequence
$\kappa$, 
\begin{equation}\label{DA_OO}
\begin{array}{l}
\displaystyle
\int_{{\mathcal W}^k(x^{(k+1)})}
\Lambda^{(k)}(x^{(k+1)},x^{(k)})\,J_\kappa(x^{(k)};\theta)\,\mathrm{d}x^{(k)}\\
\quad \displaystyle
=\,J_\kappa(x^{(k+1)};\theta)\,\frac{\Gamma((k+1)\,\theta)}{\Gamma(\theta)} 
\prod_{i=1}^k \,\frac{\Gamma((k+1-i)\,\theta+\kappa_i\big)}{\Gamma((k+2-i)\,\theta+\kappa_i)},
\end{array}
\end{equation}
where we use the convention $\kappa_i=0$ whenever $i$ exceeds the length of $\kappa$.
\end{proposition}

\section{Proof of Theorem \ref{main_thm}} \label{sec:main}

This section is devoted to the proof of Theorem \ref{main_thm}. For $\beta=2$ the result of Theorem \ref{main_thm} has been established in \cite[Section 3]{Wa} by a direct computation exploiting the explicit formulas for the transition densities of the $(\beta=2)$-Dyson Brownian motions. Later, the statement of Theorem \ref{main_thm} was shown to hold for $\beta\ge1$ in \cite[Proposition 1.3]{GS}. The proof there relies on a construction of a sequence of continuous time Markov chains realizing a discrete version of the intertwining \eqref{main_intertw} in the sense of \cite{DF}. This sequence of Markov chains is then shown to be tight, with every limit point realizing the intertwining \eqref{main_intertw}. As announced in the introduction, our proof is direct and uses only the properties of Jack symmetric polynomials stated in Section \ref{sec:Jacks}. 

\medskip

\noindent\textbf{Proof of Theorem \ref{main_thm}. Step 1.} To obtain the theorem we are going to show that for a measure-determining class of test functions, the integrals of these functions with respect to the probability measures on both sides of \eqref{main_intertw} are the same. We start by proving the corresponding identity with the Jack symmetric polynomials in $k$ variables as test functions and the semigroups $P^{(k)}(t)$, $t\ge0,$ and $P^{(k+1)}(t)$, $t\ge0,$ in \eqref{main_intertw} replaced by their respective generators. More specifically, with 
\begin{equation}
\label{DBM_gen}
{\mathcal A}^{(k)}:=\frac{1}{2}\,\sum_{i=1}^k \frac{\partial^2}{\partial z_i^2}+\theta\,\sum_{i=1}^k \, \sum_{\substack{1\le j\le k \\ j\neq i}} \frac{1}{z_i-z_j}\,\frac{\partial}{\partial z_i},
\end{equation}
our first aim is to show 
\begin{equation}\label{intertw_gen}
L^{(k)}\,{\mathcal A}^{(k)}\,J_\kappa(\,\cdot\,;\theta)
={\mathcal A}^{(k+1)}\,L^{(k+1)}\,J_\kappa(\,\cdot\,;\theta). 
\end{equation}

\smallskip

The key observation in the evaluation of ${\mathcal A}^{(k)}\,J_\kappa(\,\cdot\,;\theta)$ is that 
\begin{equation}\label{eq:commutator}
{\mathcal A}^{(k)}={\mathcal B}_1\,{\mathcal B}_2-{\mathcal B}_2\,{\mathcal B}_1
\end{equation}
with the operators ${\mathcal B}_1,{\mathcal B}_2$ of Proposition \ref{prop_B1B2}. This is indicated in \cite[comment after equation (2.13e)]{BF1} without proof and can be seen as follows: 
\begin{eqnarray*}
&& {\mathcal B}_1\,{\mathcal B}_2-{\mathcal B}_2\,{\mathcal B}_1 \\
&& =\frac{1}{2}\sum_{m,i=1}^k \!\bigg(\frac{\partial}{\partial z_m} z_i\frac{\partial^2}{\partial z_i^2} \!-\! z_i \frac{\partial^2}{\partial z_i^2} \frac{\partial}{\partial z_m}\bigg)\!+\!\theta \sum_{\substack{1\le m,i,j\le k \\ j\neq i}} \!\bigg(\frac{\partial}{\partial z_m} \frac{z_i}{z_i-z_j} \frac{\partial}{\partial z_i}\!-\!\frac{z_i}{z_i-z_j} \frac{\partial}{\partial z_i} \frac{\partial}{\partial z_m}\bigg) \\
&& =\frac{1}{2}\,\sum_{m,i=1}^k \delta_{mi}\,\frac{\partial^2}{\partial z_i^2}+\theta\,\sum_{\substack{1\le m,i,j\le k \\ j\neq i}}
\frac{(z_i-z_j)\,\delta_{mi}-z_i\,(\delta_{mi}-\delta_{mj})}{(z_i-z_j)^2}\,\frac{\partial}{\partial z_i}={\mathcal A}^{(k)},
\end{eqnarray*}
where $\delta_{mi},\delta_{mj}$ stand for Kronecker deltas.

\medskip

Combining \eqref{eq:commutator} with the identities of Proposition
\ref{prop_B1B2} we find that 
${\mathcal A}^{(k)}\,J_\kappa(x^{(k)};\theta)$ is equal to 
\begin{equation}\label{DBM_gen_eval}
\begin{array}{l} 
\displaystyle 
\,J_\kappa(1_k;\theta)\,\sum_{i=1}^l \left[\binom{\kappa}{\kappa_{(i)}}_\theta\,(\kappa_i-1+(k-i)\,\theta)\,\frac{{\mathcal B}_1\,J_{\kappa_{(i)}}(x^{(k)};\theta)}{J_{\kappa_{(i)}}(1_k;\theta)} \right.
\\
\qquad \qquad \qquad \qquad \left. \displaystyle -J_\kappa(1_k;\theta)\,\sum_{i=1}^l \binom{\kappa}{\kappa_{(i)}}_\theta\,\frac{{\mathcal B}_2\,J_{\kappa_{(i)}}(x^{(k)};\theta)}{J_{\kappa_{(i)}}(1_k;\theta)} \right]
\\
\displaystyle 
\quad = J_\kappa(1_k;\theta)\,\sum_{i=1}^l \binom{\kappa}{\kappa_{(i)}}_\theta (\kappa_i-1+(k-i)\,\theta)\,\sum_{j=1}^{l_i} \binom{\kappa_{(i)}}{(\kappa_{(i)})_{(j)}}_\theta\, \frac{J_{(\kappa_{(i)})_{(j)}}(x^{(k)};\theta)}{J_{(\kappa_{(i)})_{(j)}}(1_k;\theta)} 
\\
\displaystyle 
\qquad - J_\kappa(1_k;\theta)\,\sum_{i=1}^l \binom{\kappa}{\kappa_{(i)}}_\theta \, \sum_{j=1}^{l_i} \binom{\kappa_{(i)}}{(\kappa_{(i)})_{(j)}}_\theta ((\kappa_{(i)})_j-1+(k-j)\,\theta)\,\frac{J_{(\kappa_{(i)})_{(j)}}(x^{(k)};\theta)}{J_{(\kappa_{(i)})_{(j)}}(1_k;\theta)} 
\\
\displaystyle 
\quad =\,J_\kappa(1_k;\theta)\,\sum_{i=1}^l\,\sum_{j=1}^{l_i}
\binom{\kappa}{\kappa_{(i)}}_\theta \binom{\kappa_{(i)}}{(\kappa_{(i)})_{(j)}}_\theta (\kappa_i-(\kappa_{(i)})_j+(j-i)\,\theta)\,\frac{J_{(\kappa_{(i)})_{(j)}}(x^{(k)};\theta)}{J_{(\kappa_{(i)})_{(j)}}(1_k;\theta)},
\end{array}
\end{equation}
where $l$ and $l_i$ denote the lengths of the sequences $\kappa$ and $\kappa_{(i)}$, respectively. Now, Proposition \ref{prop_DA_OO} yields
\begin{equation}\label{LHS_gen_eval}
\begin{array}{l}
\displaystyle\int_{{\mathcal W}^k(x^{(k+1)})} \Lambda(x^{(k+1)},x^{(k)})\,{\mathcal
  A}^{(k)}\,J_\kappa(x^{(k)};\theta)\,\mathrm{d}x^{(k)} \\
\displaystyle 
\quad =J_\kappa(1_k;\theta)\,\sum_{i=1}^l\,\sum_{j=1}^{l_i} \left[ \binom{\kappa}{\kappa_{(i)}}_\theta \binom{\kappa_{(i)}}{(\kappa_{(i)})_{(j)}}_\theta 
(\kappa_i-(\kappa_{(i)})_j +(j-i)\theta)  \frac{J_{(\kappa_{(i)})_{(j)}}(x^{(k+1)};\theta)}{J_{(\kappa_{(i)})_{(j)}}(1_k;\theta)}\right. \\
\qquad \qquad \qquad \qquad \qquad \qquad \left. \displaystyle
 \cdot \frac{\Gamma(\!(k\!+\!1)\theta)}{\Gamma(\theta)}\prod_{m=1}^k
  \frac{\Gamma(\!(k+1-m)\theta\!+\!(\!(\kappa_{(i)})_{(j)})_m)}{\Gamma(\!(k+2-m)\theta\!+\!(\!(\kappa_{(i)})_{(j)})_m)}
\right], 
\end{array}
\end{equation}  
with the convention $(\!(\kappa_{(i)})_{(j)})_i=0$ whenever $i$ exceeds the length of $(\kappa_{(i)})_{(j)}$.
On the other hand, applying ${\mathcal
  A}^{(k+1)}$ to both sides of \eqref{DA_OO} and using 
\eqref{DBM_gen_eval}, with $k$ replaced by $k+1$,  we obtain 
\begin{equation}\label{RHS_gen_eval}
\begin{array}{l}
\displaystyle
{\mathcal A}^{(k+1)} \! \int_{{\mathcal W}^k(x^{(k+1)})} \!
\Lambda^{(k)}(x^{(k+1)},x^{(k)})
J_\kappa(x^{(k)};\theta)\mathrm{d}x^{(k)} \\
\displaystyle
\quad =
J_\kappa(1_{k+1};\theta) \sum_{i=1}^l \sum_{j=1}^{l_i} \!
\left[ \binom{\kappa}{\kappa_{(i)}}_{\!\!\theta} \binom{\kappa_{(i)}}{(\kappa_{(i)})_{(j)}}_{\!\!\theta} 
\cdot(\kappa_i-(\kappa_{(i)})_j+(j-i)\theta)
\frac{J_{(\kappa_{(i)})_{(j)}}(x^{(k+1)};\theta)}{J_{(\kappa_{(i)})_{(j)}}(1_{k+1};\theta)}\,
\right. 
\\
\displaystyle \left. 
\qquad \qquad \qquad \qquad \qquad \qquad \qquad \cdot
\frac{\Gamma(\!(k+1)\theta)}{\Gamma(\theta)}\prod_{m=1}^k
\frac{\Gamma(\!(k+1-m)\theta+\kappa_m\big)}{\Gamma(\!(k+2-m)\theta+\kappa_m)}\right]. 
\end{array}
\end{equation}

\smallskip

We claim that the right-hand sides of \eqref{LHS_gen_eval} and
\eqref{RHS_gen_eval} are identical. Indeed, an inspection of the
coefficients of $J_{(\kappa_{(i)})_{(j)}}(x^{(k+1)};\theta)$ in both
expressions shows that it suffices to check that 
\begin{equation}\label{Jack_pochhammer}
\begin{array}{l}
\displaystyle 
\frac{J_\kappa(1_k;\theta)}{J_\kappa(1_{k+1};\theta)}\prod_{m=1}^k
\frac{\Gamma(\!(k+2-m)\theta+\kappa_m)}{\Gamma(\!(k+1-m)\theta+\kappa_m)} \\
\qquad
\displaystyle 
=\,\frac{J_{(\kappa_{(i)})_{(j)}}(1_k;\theta)}{J_{(\kappa_{(i)})_{(j)}}(1_{k+1};\theta)} \!\prod_{m=1}^k \frac{\Gamma(\!(k+2-m)\theta+(\!(\kappa_{(i)})_{(j)})_m)}{\Gamma(\!(k+1-m)\theta+(\!(\kappa_{(i)})_{(j)})_m)}.
\end{array}
\end{equation}
However, this holds because, using the second formula in
\eqref{norm_formula}, both sides of \eqref{Jack_pochhammer} can be
shown to be equal to 
\begin{equation*}
\prod_{m=1}^k \frac{\Gamma(\!(k+2-m)\theta)}{\Gamma(\!(k+1-m)\theta)}=\frac{\Gamma(\!(k+1)\theta)}{\Gamma(\theta)}.
\end{equation*}
This finishes the proof of \eqref{intertw_gen}. 

\bigskip

\noindent\textbf{Step 2.} In this step we obtain an exponential moment
estimate for the left-hand side of \eqref{main_intertw} which will
allow us to use the moment method below. More specifically, we are
going to verify that for every $\eta>0$ the expectation
$\E^{x^{(k)}}\!\big[e^{\eta\|X^{(k)}(t)\|}\big]$ is bounded above by a
finite constant uniformly on compact sets of $(x^{(k)},t)$ in
$\mathbb{R}^k \times [0,\infty)$.  Here the notation is that of \eqref{eq:DBM} and $\|\cdot\|$ stands for the Euclidean norm on $\rr^k$.

\medskip

From It\^o's formula it follows that
\begin{eqnarray*}
\mathrm{d}\big\|X^{(k)}(t)\big\|^2 & = &  \beta
\sum_{\substack{1\le i,j\le k \\ i\neq j}} \frac{X^{(k)}_i(t)}{X^{(k)}_i(t)-X^{(k)}_j(t)}\,\mathrm{d}t + \sum_{i=1}^k 2 X^{(k)}_i(t)\,\mathrm{d}B^{(k)}_i(t) + k\,\mathrm{d}t \\
& = & \bigg(\beta\,\binom{k}{2}+k\bigg)\,\mathrm{d}t+\sum_{i=1}^k 2 X^{(k)}_i(t)\,\mathrm{d}B^{(k)}_i(t).
\end{eqnarray*}
By computing the quadratic variation of $\|X^{(k)}\|^2$ we conclude
further that $\|X^{(k)}\|^2$ is a squared Bessel process of dimension
$\beta\binom{k}{2}+k$. The desired exponential moment estimate now
follows from the Gaussian tail estimate for the latter (see, e.g.,  \cite[Section XI.1]{RY}).

\bigskip

\noindent\textbf{Step 3.} Next, we extend \eqref{intertw_gen} from the generators ${\mathcal A}^{(k)}$ and ${\mathcal A}^{(k+1)}$ to the semigroups $P^{(k)}(t)$, $t\ge0,$ and $P^{(k+1)}(t)$, $t\ge0$. To this end, we apply It\^o's formula to $J_\kappa(X^{(k)}(t);\theta)$ and take expectations on both sides of the resulting equation (keeping in mind the exponential moment estimate of Step 2) to obtain
\begin{equation}\label{Kolmogorov_Jack}
P_k(t) J_\kappa(\,\cdot\,;\theta)=J_\kappa(\,\cdot\,;\theta)+\int_0^t
P_k(s)\,{\mathcal A}^{(k)}
J_\kappa(\,\cdot\,;\theta)\,\mathrm{d}s,\quad t\ge0, 
\end{equation}
for any Jack symmetric polynomial $J_\kappa(\,\cdot\,;\theta)$. 

\medskip

Next, we recall from \eqref{DBM_gen_eval} that ${\mathcal A}^{(k)} J_\kappa(\,\cdot\,;\theta)$ is given by a finite linear combination of Jack symmetric polynomials $J_\mu(\,\cdot\,;\theta)$ satisfying $\mu_i\le\kappa_i$, $i=1,2,\ldots,l$ (with the convention $\mu_i=0$ whenever $i$ exceeds the length of $\mu$). Clearly, each $J_\mu(\,\cdot\,;\theta)$ also obeys the integral equation \eqref{Kolmogorov_Jack} and we can evaluate ${\mathcal A}^{(k)} J_\mu(\,\cdot\,;\theta)$ therein via \eqref{DBM_gen_eval}. Iterating this procedure we end up with a system of linear ordinary integral equations, whose unique solution is given by a matrix exponential. More specifically, one can represent the action of ${\mathcal A}^{(k)}$ on the finite-dimensional vector space spanned by the Jack symmetric polynomials in consideration by a matrix, whose exponential evaluated on the vector of initial values $J_\mu(0;\theta)$ yields the solution to the system of the linear ordinary integral equations. This, and the same consideration with $(k+1)$ in place of $k$, show that the identity
\begin{equation}\label{Jack_intertwining}
L^{(k)}\,P^{(k)}(t)\,J_\kappa(\,\cdot\,;\theta)
=P^{(k+1)}(t)\,L^{(k)}\,J_\kappa(\,\cdot\,;\theta),\quad t\ge0, 
\end{equation}
is the result of \eqref{intertw_gen} and the following elementary lemma.

\begin{lemma}
Suppose that three $n\times n$ matrices $M_1$, $M_2$, $M_3$ satisfy $M_1\,M_2=M_3\,M_1$. Then, $M_1\,e^{tM_2}=e^{tM_3}\,M_1$, $t\ge0$.
\end{lemma}

\smallskip

\noindent\textbf{Step 4.} To finish the proof of the theorem we recall from Definition \ref{def:Jack} that the lexicographically leading monomial with a non-zero coefficient in $J_\kappa(\,\cdot\,;\theta)$ is $z_1^{\kappa_1}z_2^{\kappa_2}\cdots z_l^{\kappa_l}$. This shows that every symmetric polynomial in $k$ variables can be written as a finite linear combination of Jack symmetric polynomials in $k$ variables and, thus, \eqref{Jack_intertwining} extends to all symmetric polynomials. Moreover, every probability measure $\gamma$ on $\mathcal{W}^{(k)}$ gives rise to a symmetrized probability measure on $\rr^k$ via
\begin{equation*}
\gamma^{\text{symm}}(\mathrm{d}z_1,\mathrm{d}z_2,\ldots,\mathrm{d}z_k):=\frac{1}{k!}\,\gamma(\mathrm{d}z_{(1)},\mathrm{d}z_{(2)},\ldots,\mathrm{d}z_{(k)}),
\end{equation*}
where $z_{(1)}\le z_{(2)}\le\cdots\le z_{(k)}$ are the order statistics of the vector $(z_1,z_2,\ldots,z_k)$. In addition, for every polynomial $p$ in $k$ variables,
\begin{equation*}
\begin{split}
\int_{\rr^k} p\,\mathrm{d}\gamma^{\text{symm}}=\int_{\rr^k} \frac{1}{k!}\,\sum_{\sigma\in S_k} p(z_{\sigma(1)},z_{\sigma(2)},\ldots,z_{\sigma(k)})\,\mathrm{d}\gamma^{\text{symm}} \\
=\int_{\mathcal{W}^k} \frac{1}{k!}\,\sum_{\sigma\in S_k} p(z_{\sigma(1)},z_{\sigma(2)},\ldots,z_{\sigma(k)})\,\mathrm{d}\gamma,
\end{split}
\end{equation*}
where $S_k$ is the set of permutations of $\{1,2,\ldots,k\}$. Using this observation for the probability measures on both sides of \eqref{main_intertw} we see that all moments of their symmetrized versions coincide. Therefore, in view of \cite[Theorem 1.1 and the remark following it]{dJ} applied to the left-hand side of \eqref{main_intertw} (recall the exponential moment estimate of Step 2), the symmetrized versions of the two sides of \eqref{main_intertw} must be equal. The theorem readily follows. \ep

\section{Random matrix proof for $\beta=1$} \label{sec:beta1}

In this section we give a much simpler proof of Theorem \ref{main_thm} for the case $\beta=1$, which relies on the random matrix interpretation of $\beta=1$ Dyson Brownian motions. The same proof applies to the cases $\beta=2,4$ as well and is omitted. 

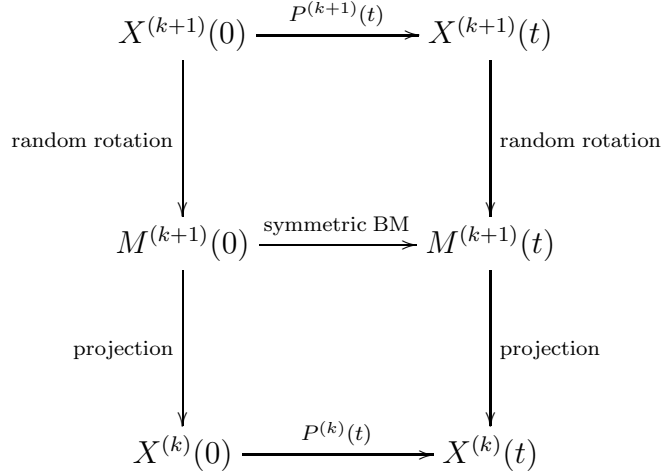
\begin{figure}[t]
\centerline{
\xymatrix@=5em{
X^{(k+1)}(0) \ar[d]_{\text{random\;rotation}} \ar[r]^{P^{(k+1)}(t)} & X^{(k+1)}(t) \ar[d]^{\text{random\;rotation}} \\
M^{(k+1)}(0) \ar[d]_{\text{projection}} \ar[r]^{\text{symmetric\;BM}} & M^{(k+1)}(t) \ar[d]^{\text{projection}} \\
X^{(k)}(0) \ar[r]^{P^{(k)}(t)} & X^{(k)}(t)
}
} 
\caption{Illustration of the proof of Theorem \ref{main_thm} for $\beta=1$.}
\label{fig:beta1}
\end{figure}

\medskip

\noindent\textbf{Proof of Theorem \ref{main_thm} for $\beta=1$.} The strategy of the proof is to construct a process realization of the intertwining \eqref{main_intertw}. To this end, we start with an arbitrary vector $X^{(k+1)}(0)\in\mathcal W^{(k+1)}$ and form the $(k+1)\times(k+1)$ diagonal matrix $D^{(k+1)}(0)$ with the components of $X^{(k+1)}(0)$ on the diagonal. Next, we sample a random $(k+1)\times(k+1)$ orthogonal matrix $O^{(k+1)}(0)$ according to the Haar measure on the orthogonal group $O(k+1)$ and introduce the symmetric matrix
\begin{equation*}
M^{(k+1)}(0):=O^{(k+1)}(0)\,D^{(k+1)}(0)\,O^{(k+1)}(0)^{-1}. 
\end{equation*}
We further define $X^{(k)}(0)\in\mathcal W^{(k)}$ as the vector of the ordered eigenvalues of the $k\times k$ top left corner of $M^{(k+1)}(0)$ (see Figure \ref{fig:beta1} for an illustration). 

\medskip

At this point, we recall that for a $(k+1)\times(k+1)$ random matrix from the Gaussian orthogonal ensemble (GOE, see e.g. \cite[Section 1.1]{Fo2} for more details) the eigenvalues are independent from the eigenvectors, with the matrix of eigenvectors being Haar distributed on the orthogonal group $O(k+1)$. Hence, the conditional probability density of $X^{(k)}(0)$ given $X^{(k+1)}(0)=x^{(k+1)}$ is the same as the conditional probability density of the ordered eigenvalues of the $k\times k$ top left corner of a $(k+1)\times(k+1)$ GOE matrix given that the vector of the ordered eigenvalues of the full matrix is $x^{(k+1)}$. That conditional probability density is known to be $\Lambda^{(k)}(x^{(k+1)},\,\cdot\,)$ (see e.g. the remark following the proof of Proposition 4.3.3 in \cite{Fo2}, as well as the proof of that proposition).  

\medskip

Next, we let the entries $(M^{(k+1)}(0))_{i,j}$, $1\le i\le j\le N,$
of $M^{k+1}(0)$ evolve according to independent standard Brownian
motions and write $M^{(k+1)}(t)$ for the completion to a symmetric
matrix of the result of such an evolution. In addition, we define
$X^{(k+1)}(t)$ and $X^{(k)}(t)$ to be the (random) vectors comprised of the ordered eigenvalues of $M^{(k+1)}(t)$ and the $k\times k$ top left corner of $M^{(k+1)}(t)$, respectively. Clearly, the law of $M^{(k+1)}(t)$ is invariant under conjugation by orthogonal matrices, so that its eigenvalues are independent from its eigenvectors, with the matrix of the latter being Haar distributed on the orthogonal group $O(k+1)$. As before, we conclude that the conditional probability density of $X^{(k)}(t)$ given $X^{(k+1)}(t)=x^{(k+1)}$ is $\Lambda^{(k)}(x^{(k+1)},\,\cdot\,)$.

\medskip

Finally, appealing to \cite[Theorem 4.3.2]{AGZ} (see also the original reference \cite[p. 123]{Mc}) we find that the evolutions of the processes $X^{(k+1)}$ and $X^{(k)}$ are governed by the semigroups $P^{(k+1)}(t)$, $t\ge0,$ and $P^{(k)}(t)$, $t\ge0$, respectively. Thus, both sides of \eqref{main_intertw} describe the conditional distribution of $X^{(k)}(t)$ given $X^{(k+1)}(0)=x^{(k+1)}$. \ep

\section{Intertwining of Dyson Ornstein-Uhlenbeck processes} \label{sec:OU}

In this last section we give the proof of Proposition \ref{prop_OU}. 

\medskip

\noindent\textbf{Proof of Proposition \ref{prop_OU}.} We follow the strategy of the proof of Theorem \ref{main_thm} and start by verifying 
\begin{equation*}
L^{(k)}\,\widetilde{\mathcal A}^{(k)}\,J_\kappa(\,\cdot\,;\theta)=\widetilde{\mathcal A}^{(k+1)}\,L^{(k)}\,J_\kappa(\,\cdot\,;\theta)
\end{equation*} 
for all Jack symmetric polynomials $J_\kappa(\,\cdot\,;\theta)$. Here 
\begin{equation*}
\widetilde{\mathcal A}^{(k)}={\mathcal A}^{(k)}-\frac{1}{2}\,{\mathcal B}_3
\end{equation*}
is the generator of the $k$-dimensional Dyson Ornstein-Uhlenbeck
process (see \eqref{DBM_gen} and \eqref{B3def} for the definitions of
${\mathcal A}^{(k)}$ and ${\mathcal B}_3$). The identities
\eqref{intertw_gen} and \eqref{B3action} imply 
\begin{equation*}
\begin{split}
L^{(k)}\,\widetilde{\mathcal A}^{(k)}\,J_\kappa(\,\cdot\,;\theta)
&=L^{(k)}\,{\mathcal A}^{(k)}\,J_\kappa(\,\cdot\,;\theta)-\frac{1}{2}\,L^{(k)}\,{\mathcal B}_3\,J_\kappa(\,\cdot\,;\theta) \\
&={\mathcal A}^{(k+1)}\,L^{(k)}\,J_\kappa(\,\cdot\,;\theta)-\frac{|\kappa|}{2}\,L^{(k)}\,J_\kappa(\,\cdot\,;\theta) 
=\widetilde{\mathcal A}^{(k+1)}\,L^{(k)}\,J_\kappa(\,\cdot\,;\theta),
\end{split}
\end{equation*}
where the last equality follows from $L^{(k)}\,J_\kappa(\,\cdot\,;\theta)$ being a multiple of $J_\kappa(\,\cdot\,;\theta)$ (see \eqref{DA_OO}) and another application of \eqref{B3action}. 

\medskip

Next, we check that for every $\eta>0$ the exponential moment $\E^{y^{(k)}}\!\big[e^{\eta\|Y^{(k)}(t)\|}\big]$ of the $k$-dimensional Dyson Ornstein-Uhlenbeck process can be bounded uniformly on compact sets of $(y^{(k)},t)$. To this end, we apply It\^o's formula to find
\begin{equation*}
\begin{split}
\mathrm{d}\big\|Y^{(k)}(t)\big\|^2 \!=\! \beta \!
\sum_{\substack{1\le i,j\le k \\ i\neq j}} \! \frac{Y^{(k)}_i(t)}{Y^{(k)}_i(t)-Y^{(k)}_j(t)}\,\mathrm{d}t \!-\! \big\|Y^{(k)}(t)\big\|^2\,\mathrm{d}t \!+\! \sum_{i=1}^k 2 Y^{(k)}_i(t)\,\mathrm{d}B^{(k)}_i(t) \!+\! k\,\mathrm{d}t \\
=\bigg(\beta\,\binom{k}{2}+k-\big\|Y^{(k)}(t)\big\|^2\bigg)\,\mathrm{d}t+\sum_{i=1}^k 2 Y^{(k)}_i(t)\,\mathrm{d}B^{(k)}_i(t).
\end{split}
\end{equation*}
In other words, $\|Y^{(k)}\|^2$ solves the stochastic differential equation 
\begin{equation*}
\mathrm{d}R(t)=\bigg(\beta\,\binom{k}{2}+k-R(t)\bigg)\,\mathrm{d}t+2\sqrt{R(t)}\,\mathrm{d}W(t),
\end{equation*}
where $W$ is a standard Brownian motion. Thus, by \cite[Proposition 5.2.18]{KS} the process $\|Y^{(k)}\|^2$ can be dominated pathwise by a squared Bessel process of dimension $\beta\binom{k}{2}+k$. The desired exponential moment estimate readily follows from the Gaussian tail estimate for the latter (see e.g. \cite[Section XI.1]{RY}).

\medskip

To conclude the proof of the proposition it remains to repeat Steps 3
and 4 of the proof of Theorem \ref{main_thm} word-by-word. 
We omit the details. \ep

\bigskip\bigskip\bigskip


\end{document}